\documentclass[12pt,letterpaper,twoside]{article}
\usepackage[centertags]{amsmath}
\usepackage{graphicx,indentfirst,amsmath,amsfonts,amssymb,amsthm,newlfont}
\font\Goth=yinitas scaled \magstep0

\newcommand{\Gth}[1]{\lower2mm\hbox{\Goth #1}}

\usepackage[latin1]{inputenc}
\def\al{\alpha}

\def\be{\beta}
\def\De{\Delta}

\def\l1{{\lambda}_1}

\newcommand{\f}{\frac}
\newcommand{\lh}{\Delta_{H^{1}}}

\def\x1{{\xi }_{xx}}
\def\x2{{\xi }_{yy}}
\def\x3{{\xi }_{xy}}

\def\e1{{\eta }_{xx}}
\def\e2{{\eta }_{yy}}
\def\e3{{\eta }_{xy}}

\def\kd{\partial}

\newcommand{\ds}{\displaystyle }
\newtheorem{theorem}{Theorem}
\newtheorem{lemma}{\bf Lemma}

\newcommand{\beqn}{\begin{eqnarray*}}
\newcommand{\eeqn}{\end{eqnarray*}}
\newcommand{\beqnn}{\begin{eqnarray}}
\newcommand{\eeqnn}{\end{eqnarray}}
\newcommand{\p}{\partial}
\newcommand{\bb}{\begin{equation}}
\newcommand{\ee}{\end{equation}}
\newcommand{\ba}{\begin{array}}
\newcommand{\ea}{\end{array}}
\newcommand{\R}{\mathbb{R}}
\newcommand{\N}{\mathbb{N}}

\newcommand{\X}{\tilde{X}}
\newcommand{\Y}{\tilde{Y}}
\newcommand{\W}{W_\be}
\begin{document}
\pagenumbering{arabic}
\title{\huge \bf Symmetry Coefficients of Semilinear Partial Differential Equations}
\author{\rm \large Igor Leite Freire and Antonio Carlos Gilli Martins\\
\\
\it Instituto de Matem\'atica,
Estat\'\i stica e \\ \it Computa\c c\~ao Cient\'\i fica - IMECC \\
\it Universidade Estadual de Campinas - UNICAMP \\ \it C.P.
$6065$, $13083$-$970$ - Campinas - SP, Brasil
\\ \rm E-mail: igor@ime.unicamp.br\\ \ \ \ \ gilli@ime.unicamp.br }
\date{\ }
\maketitle
\vspace{1cm}
\begin{abstract}
We show that for any semilinear partial differential equation of order $m$, the infinitesi-mals of the independent variables depend only on the independent variables and, if $m>1$ and the equation is also linear in its derivatives of order $m-1$ of the dependent variable, then the infinitesimal of the dependent variable is at most linear on the dependent variable. Many examples of important partial differential equations in Analysis, Geometry and Mathematical - Physics are given in order to enlighten the main result.
\end{abstract}
\vskip 1cm
\begin{center}
{2000 AMS Mathematics Classification numbers:\vspace{0.2cm}\\
35H10, 58J70\vspace{0.2cm} }
\end{center}
\pagenumbering{arabic}
\newpage

\section{Introduction}

Let $x\in M\subseteq\mathbb{R}^{n}$, $u:M\rightarrow\mathbb{R}$ a smooth
function and $k\in\mathbb{N}$. We use $\partial^{k} u$ to denote the jet bundle corresponding to all $k$th partial derivatives of $u$ with respect to $x$. We simply denote $\p^{1}u$ by $\p u$. A partial differential equation (PDE) of order $m$ is a relation $F(x,u,\p u,\cdots,\p^{m}u)=0$.

If there exists an operator
\bb\label{linop}L_{m}:=A^{i_{1}\cdots i_{m}}(x)\f{\p^{m}}{\p x^{i_{1}}\cdots \p x^{i_{m}}}\ee
and a relation $f(x,u,\p u,\cdots,\p^{m-1}u)$ such that $F=Lu+f(x,u,\p u,\cdots,\p^{m-1}u)$, then $F=0$ is said to be a \textit{semilinear partial differential equation} (SPDE). In this article we use the Einstein summation convention.

Partial differential equations are used to model many different kinds of phenomena in science and engeneering. Linear equations give mathematical description for physical, chemical or biological processes in a first approximation only. In order to have a more detailed and precise description a mathematical model needs to incoporate nonlinear terms. Nonlinear equations are difficult to solve analytically. However, in the end of century $XIX$ Sophus Lie developed a method that is widely useful to obtain solutions of a differential equation. This method is currently called \textit{Lie point symmetry theory}. Some applications of this method in (nonlinear) differential equations can be found in \cite{baks, pbh, afby,bl, yb2, gc, ds, cl, yii,ilf, ib, la, smpg, ol, sv}.

Lie used group properties of differential equations in order to actually solve them, i.e., to construct their exact solutions. Nowadays symmetry reductions are one of the most powerful tools for solving nonlinear PDEs.

A Lie point symmetry\footnote{In fact, a Lie point symmetry is given by the exponential map $(\exp{S})(x,u)=:(x^{\ast},u^{\ast})\in\R^{n}\times\R$. We are identifying the point transformation with its generator.} of a PDE of order
$m$ is a vector field
\bb\label{simetria}
S=\xi^{i}(x,u)\frac{\partial}{\partial x^{i}}+\eta(x,u)\frac{\partial
}{\partial u}%
\ee
on $M\times\mathbb{R}$ such that $S^{(m)}F=0$ when $F=0$ and
\[
S^{(m)}:=S+\eta^{(1)}_{i}(x,u,\partial u)\frac{\partial}{\partial u_{i}}%
+\cdots+\eta^{(m)}_{i_{1}\cdots i_{m}}(x,u,\partial u,\cdots,\partial^{m}
u)\frac{\partial}{\partial u_{i_{1}\cdots i_{m}}}%
\]
is the extended symmetry on the jet space $(x,u,\partial u,\cdots,
\partial^{k} u)$.

The functions $\eta^{(j)}(x,u,\partial u,\cdots,\partial^{j} u)$, $1\leq j\leq
m$, are given by
\bb\label{etas}
\begin{array}
[c]{l c l}%
\eta^{(1)}_{ i} &: = & D_{i}\eta-(D_{i}\xi^{j})u_{j},\\
&  & \\
\eta^{(j)}_{i_{1}\cdots i_{j}} &: = & D_{i_{j}}\eta^{(j-1)}_{i_{1}\cdots
i_{j-1}}-(D_{i_{j}}\xi^{l})u_{i_{1}\cdots i_{j-1}l},\;2\leq j\leq m,
\end{array}
\ee
where
$$D_{i}:=\f{\p}{\p x^{i}}+u_{i}\f{\p}{\p u}+u_{ij}\f{\p}{\p u_{j}}+\cdots+u_{ii_{1}\cdots i_{m}}\f{\p}{\p u_{i_{1}\cdots i_{m}}}+\cdots $$
is the \textit{total derivative operator}. We shall not present more preliminaries concerning the Lie point symmetries of
differential equations supposing that the reader is familiar with the basic notions and methods of group analysis \cite{bl,ib,ol}.

In \cite{b1}, Bluman proved some relations between symmetry coefficients wich simplify drasticaly their calculus. 
He worked with a PDE of the form
\bb\label{nonli}A^{i_{1}\cdots i_{m}}(x,u)u_{i_{1}\cdots i_{m}}+f(x,u,\p u,\cdots,\p^{m-1}u)=0.\ee

Depending on the relations between coefficients $A^{i_{1}\cdots i_{m}}(x,u)$ of equation (\ref{nonli}), Bluman's theorems gives us conditions to determine, \textit{a priori}, if the coefficient $\xi^{i}$ depends or not of $u$, and in many cases, it also gives us some information about the dependence of coefficient $\eta$ with respect to $u$.

The Bluman's theorems can be used in a more general results where it is necessary to find the infinitesimals $\xi^{i}$ and $\eta$ of symmetry (\ref{simetria}) of a quaselinear equation (\ref{nonli}). However, in the special cases where the coefficients $A^{i_{1}\cdots i_{m}}(x,u)$ do not depend of the dependent variable $u$, the proof that the coefficients $\xi^{i}$ do not depend of $u$ and that $\eta$ is linear in $u$ is simpler than that presented by Bluman in \cite{b1}. And in this case, the equation (\ref{nonli}) becames a semilinear partial differential equation. Since many of the most important equations from Analysis, Geometry and Mathematical-Physics are SPDE, we intend to enlighten the Bluman's proof of his theorem for the semilinear case.

The purpose of this article is twofold. First, we intend to give a detailed proof of a theorem (Theorem \ref{ig}) wich gives us conditions to state the coefficients $\xi^{i}$ with respect to $u$ of a symmetry of a SPDE and, in many cases, we can conclude that $\eta$ is a linear function with respect to $u$ (see \cite{b1,bl}). 

The second purpose is to present and summarize some important PDEs arising from Analysis, Geometry and Mathematical - Physics, which are linear PDEs or SPDEs (see Section \ref{examples}), illustrating Theorem \ref{ig}. 

Our main purpose is to prove the following result:
\begin{theorem}\label{ig}
Consider the SPDE
\bb\label{spde}
A^{i_{1}\cdots i_{m}}(x)u_{i_{1}\cdots i_{m}}+f(x,u,\p u,\cdots,\p^{m-1}u)=0,
\ee
where $A^{i_{1}\cdots i_{m}}(x)$ is symmetric with respect its indeces. Suppose that the vector field $S$ given in $(\ref{simetria})$ is a symmetry of $(\ref{spde})$. Then $\xi^{i}_{u}=0$, $1\leq i\leq n$.

If $m>1$ and $f(x,u,\p u,\cdots,\p^{m-1}u)=a^{i_{1}\cdots i_{m-1}}(x)u_{i_{1}\cdots i_{m-1}}+h(x,u,\p u,\cdots,\p ^{m-2}u)$, for some function $h$, then $\eta_{uu}=0$. 
\end{theorem}

The paper is organized as follows. In section 2 we prove Theorem \ref{ig}.
In section 3 we give some examples, from Analysis, Geometry and Mathematical - Physics, illustrating the Theorem.\newpage

\section{Proof of the main results}
\

In this section, we shall prove Theorem \ref{ig}. We shall do this in three steps: first, we prove Theorem \ref{ig} when $m=1$. In this case we, at most, can conclude $\xi^{i}=\xi^{i}(x)$. The case $m=2$ is done because many of the most important equations in Analysis, Geometry and Mathematical - Physics are second order SPDE and this proof is a good way to understand the proof of arbitrary $m$, which is the third step.


\subsection{The case $m=1$}
\begin{proof}
Let $L:=A^{i}(x)\f{\p}{\p x^{i}}$ a linear operator and $f(x,u)$ a smooth function. Consider the first order semilinear partial differential equation
\bb\label{lineq}F(x,u,\p u):=Lu+f(x,u)=0.\ee
Suppose (\ref{lineq}) admits a symmetry $S$ given by (\ref{simetria}). Its first order extension is
$$S^{(1)}=\xi^{j}(x,u)\f{\p}{\p x^{j}}+\eta(x,u)\f{\p}{\p u}+(\eta_{i}(x,u)+\eta_{u}(x,u)u_{i}-\xi^{j}_{i}(x,u)u_{j}-\xi^{j}_{u}(x,u)u_{j}u_{i})\f{\p}{\p u_{i}}.$$
Apllying $S^{(1)}$ to (\ref{lineq}), we have
$$S^{(1)}F=(\xi^{i}f_{i}+\eta f_{u}+A^{i}\eta_{i})+(\xi^{i}A^{j}_{i}+A^{j}\eta_{u}-A^{i}\xi_{i}^{j})u_{j}-A^{i}\xi^{i}_{u}u_{i}u_{j}.$$

Then, by the symmetry condition (see Ibragimov \cite{ib} or Olver \cite{ol})
$$S^{(1)}F=\lambda(x,u)F$$
and since $F$ is a linear function with respect to $\p u$, we conclude that
\bb\label{eq1}A^{i}\xi^{j}_{u}u_{i}u_{j}=0\ee
Choising $i_{0}$ such that $A^{i_{0}}\neq 0$, the equation (\ref{eq1}) implies that necessarily we must have $\xi_{u}^{j}=0$. Thus, $\xi^{j}=\xi^{j}(x)$ and this conclude the proof for the case $m=1$.


\subsection{The case $m=2$}
\

Let $$F:=A^{ij}(x)u_{ij}+f(x,u,\p u)=0$$ be a SPDE and $S^{(2)}$ the second order extension of symmetry (\ref{simetria}). Then,  
$$S^{(2)}=\xi^{k}(x,u)\f{\p}{\p x^{k}}+\eta(x,u)\f{\p}{\p u}+\eta_{k}^{(1)}(x,u,\p u)\f{\p}{\p u_{k}}+\eta_{kl}^{(2)}(x,u,\p u,\p^{2} u)\f{\p}{\p u_{kl}},$$
and the coefficients in the jet spaces are given by (see equation (\ref{etas}))
$$
\eta_{k}^{(1)} =  \ds{\eta_{k}-\xi^{j}_{k}u_{j}-\xi^{j}_{u}u_{j}u_{k}+\eta_{u}u_{k}},
$$

\bb\label{coefe}\ba{l c l}
\eta_{kl}^{(2)} & = &\ds{ \eta_{kl}+\eta_{lu}u_{k}-\xi^{j}_{kl}u_{j}+\eta_{ku}u_{l}-\xi^{j}_{lu}u_{j}u_{k}-\xi^{j}_{ku}u_{j}u_{l}-\xi^{j}_{uu}u_{j}u_{k}u_{l}-\xi^{j}_{k}u_{lj}-\xi^{j}_{l}u_{jk}}\\
\\
 & & \ds{-\xi^{j}_{u}u_{j}u_{lk}-\xi^{j}_{u}u_{k}u_{jl}-\xi^{j}_{u}u_{l}u_{kj}+\eta_{u}u_{kl}+\eta_{uu}u_{k}u_{l}}
\ea
\ee

Let $F_{k}:=\frac{\p F}{\p x^{k}}$, then
$$
\ba {l c l}
S^{(2)}F & = & \ds{\xi^{k}F_{k}+\eta F_{u}+\eta_{k}F_{u_{k}}+A^{kl}\eta_{kl}+(u_{k}\eta_{u}-\xi^{j}_{k}u_{j})F_{u_{k}}+A^{kl}\eta_{kl}u_{k}-A^{kl}\xi^{j}_{kl}u_{j}}\\
\\
& & 
\ds{+A^{kl}\eta_{ku}u_{l}-A^{kl}\xi^{j}_{lu}u_{j}u_{k}+A^{kl}\eta_{uu}u_{k}u_{l}-A^{kl}\xi^{j}_{ku}u_{j}u_{l}-A^{kl}\xi^{j}_{uu}u_{j}u_{k}u_{l}}\\
\\
& & \ds{+A^{kl}\eta_{u}u_{kl}-A^{kl}\xi^{j}_{k}u_{lj}-A^{kl}\xi^{j}_{l}u_{kj}-A^{kl}\xi^{j}_{u}u_{lj}u_{k}-A^{kl}\xi^{j}_{u}u_{lk}u_{j}-A^{kl}\xi^{j}_{u}u_{l}u_{jk}}.
\ea
$$
The symmetry condition is

\bb\label{condsim}S^{(2)}F=\lambda(x,u)F.\ee

Since $F$ is a linear function in the second order derivatives of $u$, the symmetry condition (\ref{condsim}) implies that terms $u_{i}u_{jk}$ have to be zero. Then,
\bb\label{del}A^{kl}(\xi^{j}_{u}u_{lj}u_{k}+\xi^{j}_{u}u_{lk}u_{j}+\xi^{j}_{u}u_{l}u_{jk})=0.\ee
Since
$$
\ba{l c l c l c l}
u_{k}  & = & \delta^{p}_{k}u_{p}, & & u_{lj} & = & \delta^{r}_{l}\delta^{s}_{j}u_{rs},\\
\\
u_{j}  & = & \delta^{p}_{j}u_{p}, & & u_{lk} & = & \delta^{r}_{l}\delta^{s}_{k}u_{rs},\\
\\
u_{l}  & = & \delta^{p}_{l}u_{p}, & & u_{lj} & = & \delta^{r}_{j}\delta^{s}_{k}u_{rs},
\ea
$$ 
and substituting this into (\ref{del}), we have the following relation
$$(A^{kl}\xi^{j}_{u}\delta^{p}_{k}\delta^{r}_{l}\delta^{s}_{j}+A^{kl}\xi^{j}_{u}\delta^{p}_{j}\delta^{r}_{l}\delta^{s}_{k}+A^{kl}\xi^{j}_{u}\delta^{p}_{l}\delta^{r}_{j}\delta^{s}_{k})u_{p}u_{rs}=0.$$
Since the set $\{u_{j}u_{kl}\}$ is linearly independet set, the following identity must be satisfied:
\bb\label{xiu}A^{kl}\xi^{j}_{u}\delta^{p}_{k}\delta^{r}_{l}\delta^{s}_{j}+A^{kl}\xi^{j}_{u}\delta^{p}_{j}\delta^{r}_{l}\delta^{s}_{k}+A^{kl}\xi^{j}_{u}\delta^{p}_{l}\delta^{r}_{j}\delta^{s}_{k}=0.\ee
Taking $p=r=s$, we conclude that 
$$A^{pp}\xi^{p}_{u}=0.$$

Let $N_{1}$ e $N_{2}$ be the set of indeces such that $A^{pp}\neq 0$ and $A^{pp}=0$, respectively. Then, for all $i\in N_{1}$, $\xi^{i}_{u}=0$ and hence, $\xi^{i}=\xi^{i}(x)$. 

Suppose $N_{2}\neq\emptyset$. Thus, there exists $n_{0}\in N_{2}$ such that $A^{n_{0}p}\neq 0$, for some $p$. Taking $k\neq p,\; j\neq p$ and choising $s=n_{0}$ in (\ref{xiu}), we obtain
$$A^{n_{0}p}\xi^{r}_{u}=0.$$
Thus we conclude that $\xi^{r}_{u}=0$ for all $r$. 

Now, suppose that $f=b^{j}(x)u_{j}+h(x,u)$. Since $\xi^{i}_{u}=0$, we can write
$$
\ba {l c l}
S^{(2)}F & = & \ds{\xi^{k}F_{k}+\eta F_{u}+\eta_{k}F_{u_{k}}+A^{kl}\eta_{kl}+(u_{k}\eta_{u}-\xi^{j}_{k}u_{j})F_{u_{k}}+A^{kl}\eta_{kl}u_{k}-A^{kl}\xi^{j}_{kl}u_{j}}\\
\\
& & 
\ds{+A^{kl}\eta_{ku}u_{l}+A^{kl}\eta_{u}u_{kl}-A^{kl}\xi^{j}_{k}u_{lj}-A^{kl}\xi^{j}_{l}u_{kj}+A^{kl}\eta_{uu}u_{k}u_{l}}
\ea
$$
Since the SPDE $F=0$ is linear in the first and the second order derivatives of $u$, the condition (\ref{condsim}) implies that the coefficients of $u_{k}u_{l}$ have to be zero. Then, $A^{kl}\eta_{uu}=0$ and, finally, $\eta_{uu}=0$.

\subsection{The case $m>2$}
\

\begin{lemma}\label{lema}
Let $k\geq 2$. Then, there exists a function\footnote{This function is a polinomial function in $\p u, \cdots, \p^{k}u$ (see \cite{b1,bl}).} h depending of $x,u,\p u, \cdots, \p^{k}u$, such that
$$
\ba{l c l}
\eta^{(k)}_{i_{1}\cdots i_{k}}& = & h(x,u,\p u, \cdots, \p^{k} u)-\xi^{j}_{u}u_{j}u_{i_{1}\cdots i_{k}}-\xi^{j}_{u}u_{i_{1}}u_{ji_{2}-\cdots i_{k}}-\cdots\xi_{u}^{j}u_{i_{k}}u_{i_{1}\cdots i_{k-1}j}\\
\\
& &+\eta_{u}u_{i_{1}\cdots i_{k}} +\eta_{uu}(u_{i_{1}}u_{i_{2}\cdots i_{k-1}}+u_{i_{2}}u_{i_{1}i_{3}\cdots i_{k-1}}+\cdots+u_{i_{k}}u_{i_{1}\cdots i_{k-1}}).
\ea$$

\end{lemma}

\begin{proof}
We shall prove that the Lemma is valid for all $k>1$. If $k=2$, we turn back to equation (\ref{coefe}). Suppose that the result is valid to $k,\;k>2$. Then,

$$
\ba{l c l}
\eta^{(k)}_{i_{1}\cdots i_{k}}& = & h(x,u,\p u, \cdots, \p^{k} u)-\xi^{j}_{u}u_{j}u_{i_{1}\cdots i_{k}}-\xi^{j}_{u}u_{i_{1}}u_{ji_{2}\cdots i_{k}}-\cdots-\xi_{u}^{j}u_{i_{k}}u_{i_{1}\cdots i_{k-1}j}\\
\\
& & +\eta_{u}u_{i_{1}\cdots i_{k-1}}+\eta_{uu}(u_{i_{1}}u_{i_{2}\cdots i_{k-1}}+u_{i_{2}}u_{i_{1}i_{3}\cdots i_{k-1}}+\cdots+u_{i_{k}}u_{i_{1}\cdots i_{k-1}}).
\ea$$

From equation (\ref{etas}), we have, after a straigthforward calculation,
\begin{equation}\label{ind}
\ba{l c l}
\eta_{i_{1}\cdots i_{k}i_{k+1}}^{(k+1)}& = &(D_{i_{k+1}}h)-(D_{i_{k+1}}\xi^{j}_{u})u_{j}u_{i_{1}\cdots i_{k}}-(D_{i_{k+1}}\xi^{j}_{u})u_{i_{1}}u_{j i_{2}\cdots i_{k}}-\cdots\\
\\
& & -(D_{i_{k+1}}\xi^{j}_{u})u_{i_{k}}u_{i_{1}\cdots i_{k-1}j}-\xi^{j}_{i_{k+1}}u_{i_{1}\cdots i_{k}j}-\xi^{j}_{u}u_{ji_{k+1}}u_{i_{1}\cdots i_{k}}\\
\\
& &
 -\xi^{j}_{u}u_{i_{1}i_{k+1}}u_{ji_{2}\cdots i_{k}}-\xi^{j}_{u}u_{i_{2}i_{k+1}}u_{i_{1}ji_{3}\cdots i_{k}}-\cdots-\xi^{j}_{u}u_{i_{k}i_{k+1}}u_{i_{1}\cdots i_{k-1}j}\\
 \\
& &-\eta_{i_{k+1}u}u_{i_{1}\cdots i_{k}}+\eta_{uu}(u_{i_{1}i_{k+1}}u_{i_{2}\cdots i_{k}}+\cdots+u_{i_{k}i_{k+1}}u_{i_{2}\cdots i_{k}})\\
\\
& &-\xi^{j}_{u}u_{j}u_{i_{1}\cdots i_{k}i_{k+1}}-\xi^{j}_{u}u_{i_{1}}u_{ji_{2}\cdots i_{k}i_{k+1}}-\cdots -\xi^{j}_{u}u_{i_{k}}u_{i_{1}\cdots i_{k-1}ji_{k+1}}\\
\\
& & 
-\xi^{j}_{u}u_{i_{k+1}}u_{i_{1}\cdots i_{k}j}+\eta_{u}u_{i_{1}\cdots i_{k}i_{k+1}}+\eta_{uu}(u_{i_{1}}u_{i_{2}\cdots i_{k}i_{k+1}}+\cdots+u_{i_{k+1}}u_{i_{2}\cdots i_{k}})
\ea
\end{equation}
Let 
$$
\ba{l c l}
\tilde{h}(x,u,\p u,\cdots \p^{k+1}u)& := &(D_{i_{k+1}}h)-(D_{i_{k+1}}\xi^{j}_{u})u_{j}u_{i_{1}\cdots i_{k}}-(D_{i_{k+1}}\xi^{j}_{u})u_{i_{1}}u_{j i_{2}\cdots i_{k}}-\cdots\\
\\
& & -(D_{i_{k+1}}\xi^{j}_{u})u_{i_{k}}u_{i_{1}\cdots i_{k-1}j}-\xi^{j}_{i_{k+1}}u_{i_{1}\cdots i_{k}j}-\xi^{j}_{u}u_{ji_{k+1}}u_{i_{1}\cdots i_{k}}\\
\\
& &
 -\xi^{j}_{u}u_{i_{1}i_{k+1}}u_{ji_{2}\cdots i_{k}}-\cdots-\xi^{j}_{u}u_{i_{k}i_{k+1}}u_{i_{1}\cdots i_{k-1}j}-\eta_{i_{k+1}u}u_{i_{1}\cdots i_{k}}\\
 \\
 & & +\eta_{uu}(u_{i_{1}i_{k+1}}u_{i_{2}\cdots i_{k}}+\cdots+u_{i_{k}i_{k+1}}u_{i_{2}\cdots i_{k}}).
\ea
$$Then, we conclude that equation (\ref{ind}) can be written as
$$
\ba{l c l}
\eta_{i_{1}\cdots i_{k}i_{k+1}}^{(k+1)}& = & \tilde{h}-\xi^{j}_{u}u_{j}u_{i_{1}\cdots i_{k}i_{k+1}}-\xi^{j}_{u}u_{i_{1}}u_{ji_{2}\cdots i_{k}i_{k+1}}-\cdots-\xi^{j}_{u}u_{i_{k}}u_{i_{1}\cdots i_{k-1}ji_{k+1}}\\
\\
& & -\xi^{j}_{u}u_{i_{k+1}}u_{i_{1}\cdots i_{k}j}+\eta_{u}u_{i_{1}\cdots i_{k}i_{k+1}}+\eta_{uu}(u_{i_{1}}u_{i_{2}\cdots i_{k}i_{k+1}}+\cdots+u_{i_{k+1}}u_{i_{2}\cdots i_{k}}),
\ea
$$
proving the Lemma.
\end{proof}
Now, we are in position to prove the general case:
Let $$F:=A^{i_{1}\cdots i_{m}}(x)u_{i_{1}\cdots i_{m}}+f(x,u,\p u,\cdots,\p^{m-1}u)$$ and $S^{(m)}$ the extended symmetry of (\ref{simetria}). Then, by Lemma \ref{lema}, we have
\bb\label{simeta}
\ba{l c l}
S^{(m)}F & = & \xi^{j}A_{j}^{i_{1}\cdots i_{m}}u_{i_{1}\cdots i_{m}}+\xi^{j}f_{j}+\eta f_{u}+\eta^{(1)}_{i}f_{u_{i}}+\cdots \eta^{(m)}_{i_{1}\cdots i_{m-1}}f_{u_{i}\cdots i_{m-1}}+\\
\\
& & +A^{i_{1}\cdots i_{m}}[\tilde{h}-\xi^{j}_{u}u_{j}u_{i_{1}\cdots i_{k}i_{k+1}}-\xi^{j}_{u}u_{i_{1}}u_{ji_{2}\cdots i_{k}i_{k+1}}-\cdots-\xi^{j}_{u}u_{i_{k}}u_{i_{1}\cdots i_{k-1}ji_{k+1}}\\
\\
& & -\xi^{j}_{u}u_{i_{k+1}}u_{i_{1}\cdots i_{k}j}+\eta_{u}u_{i_{1}\cdots i_{k}i_{k+1}}+\eta_{uu}(u_{i_{1}}u_{i_{2}\cdots i_{k}i_{k+1}}+\cdots+u_{i_{k+1}}u_{i_{2}\cdots i_{k}})].
\ea
\ee
By the symmetry condition $S^{(m)}F=\lambda(x,u)F$, necessarily we have to have
\bb\label{coef}
A^{i_{1}\cdots i_{m}}\xi^{j}_{u}(u_{j}u_{i_{1}\cdots i_{m}}+u_{i_{1}}u_{ji_{2}\cdots i_{m}}+\cdots u_{i_{m}}u_{i_{1}\cdots i_{m-1}j})=0.
\ee
Since
$$
\ba{l c l}
u_{j}u_{i_{1}i_{2}\cdots i_{m}} & = & u_{p}u_{l_{1}l_{2}\cdots l_{m}}\delta^{pl_{1}l_{2}\cdots l_{m-1}l_{m}}_{ji_{1}i_{2}\cdots i_{m-1}i_{m}},\\
\\
u_{i_{1}}u_{ji_{2}\cdots i_{m}} & = & u_{p}u_{l_{1}l_{2}\cdots l_{m}}\delta^{pl_{1}l_{2}\cdots l_{m-1}l_{m}}_{i_{1}j\cdots i_{m-1}i_{m}},\\
\\
\vdots & \vdots & \vdots\\
\\
u_{i_{m}}u_{i_{1}i_{2}\cdots i_{m-1}j} & = & u_{p}u_{l_{1}l_{2}\cdots l_{m}}\delta^{pl_{1}l_{2}\cdots l_{m-1}l_{m}}_{i_{m}i_{1}i_{2}\cdots i_{m-1}j},
\ea
$$
where 
$$\delta^{l_{1}l_{2}\cdots l_{m-1}l_{m}}_{k_{1}k_{2}\cdots k_{m-1}k_{m}}:=\delta^{l_{1}}_{k_{1}}\delta^{l_{2}}_{k_{2}}\cdots\delta^{l_{m}}_{k_{m}}.$$

Equation (\ref{coef}) becames
\bb\label{coefgen}
A^{i_{1}\cdots i_{m}}\xi^{j}_{u}(\delta^{pl_{1}l_{2}\cdots l_{m-1}l_{m}}_{ji_{1}i_{2}\cdots i_{m-1}i_{m}}+\delta^{pl_{1}l_{2}\cdots l_{m-1}l_{m}}_{i_{1}j\cdots i_{m-1}i_{m}}+\delta^{pl_{1}l_{2}\cdots l_{m-1}l_{m}}_{i_{m}i_{1}i_{2}\cdots i_{m-1}j})u_{p}u_{l_{1}\cdots l_{m}}=0,
\ee
$p\in\{i_{1},\cdots,i_{m}\},\;i_{s}\in\{1,\cdots,n\},\;s\in\{1,\cdots,m\},\;j\in\{1,\cdots,n\}$.

Whereas the set $\{u_{p}u_{l_{1}\cdots l_{m}}\}$ is a linearly independent set, in order that equation (\ref{coefgen}) be true, we necessarily have 
$$A^{i_{1}\cdots i_{m}}\xi^{j}_{u}(\delta^{pl_{1}l_{2}\cdots l_{m-1}l_{m}}_{ji_{1}i_{2}\cdots i_{m-1}i_{m}}+\delta^{pl_{1}l_{2}\cdots l_{m-1}l_{m}}_{i_{1}j\cdots i_{m-1}i_{m}}+\delta^{pl_{1}l_{2}\cdots l_{m-1}l_{m}}_{i_{m}i_{1}i_{2}\cdots i_{m-1}j})=0.$$

Taking $l_{k}=i_{k},\;1\leq k \leq m$, such that $A^{m}:=A^{i_{1}\cdots i_{m-1}i_{m}}\neq 0$, like in the case $m=2$, we obtain
$$A^{m}\xi^{j}_{u}(\delta^{p}_{j}+\delta^{p}_{i_{1}}\delta^{l_{1}}_{j}\cdots+\delta^{p}_{i_{m}}\delta^{l_{m}}_{j})=0.$$

Since the term $\delta^{p}_{j}+\delta^{p}_{i_{1}}\delta^{l_{1}}_{j}\cdots+\delta^{p}_{i_{m}}\delta^{l_{m}}_{j}$ cannot be zero, we necessarily have $\xi^{j}_{u}=0$. Thus $\xi^{i}=\xi^{j}(x)$.



Suppose now that $f(x,u,\p u, \cdots, \p^{k-1} u)$ is linear in $\p^{k-1} u$. Then, from equation (\ref{simeta}) and the symmetry condiction (\ref{condsim}), we can see that the term 
$$A^{i_{1}\cdots i_{m}}\eta_{uu}(u_{i_{1}}u_{i_{2}\cdots i_{k}i_{k+1}}+\cdots+u_{i_{k+1}}u_{i_{2}\cdots i_{k}})=0.$$
In the same way, as in the case $m=2$, we easily conclude that $\eta_{uu}=0$. Then, there exist functions $\alpha=\alpha(x)$ and $\be=\be(x)$ such that $\eta=\al(x)u+\be(x)$.
\end{proof}

\section{Some Examples}\label{examples}
\

The following examples illustrate the Theorem 1. For another examples of equations where the theorem could be applied in order to obtain the symmetries coefficients, see \cite{baks,pbh,afby,smpg}.
\subsection{Poisson Equation}

The Poisson Equation in $\R^{n}$ is 
\bb\label{poisson}
\Delta u+f(u)=0,
\ee
where
$$\Delta:=\delta^{ij}\f{\p}{\p x^{i}}\f{\p}{\p x^{j}}$$
denotes the Laplace operator in $\R^{n}$.

When $n=2$, the group classification was obtained for Sophus Lie in the end of $XIX$ century. He proved the following result:

The widest Lie point symmetry group admitted by $(\ref{poisson})$, with
 arbitrary $f(u)$, is determined by translations
 \bb\label{e2}
 Y_1=\frac{\kd }{\kd x},\;\;\;Y_2=\frac{\kd }{\kd y}
 \ee
 and the rotation
 \bb\label{e3}
 Y_3=y\; \frac{\kd }{\kd x}-x\; \frac{\kd }{\kd y}.
 \ee

 For some special choices of $f(u)$ it can be
 expanded by operators additional to $(\ref{e2})$ and $(\ref{e3})$, which are listed
 below.
\begin{itemize}
\item If $f(u)=0$, then
 \bb\label{e4} Y_{(\xi^{1},\xi^{2})}={\xi }^1\frac{\kd }{\kd x} +{\xi }^2\frac{\kd }{\kd
 y}, \ee
 \bb\label{e5} Y_4=u\frac{\kd }{\kd u}, \;\;\; Y_{\be }= \be (x,y)
 \frac{\kd }{\kd u}, \text{ with } \De \be =0, \ee
 where ${\xi }^1={\xi }^1(x,y) $, ${\xi }^2={\xi }^2(x,y) $
 satisfy the Cauchy-Riemann equations:
 \bb\label{e8} {\xi }^1_x= {\xi }^2_y,\;\;\;{\xi }^1_y = -{\xi
 }^2_x. \ee

\item The case $f(u)=const$ can be easily reduced to the preceding
 one.

\item If $f(u)=k u$, $k\neq 0$ is a constant, we have $Y_4$ and
 $Y_{\be }$, where $\De \be +k\be=0$.

\item For $f(u)=ku^p$, $p\neq 0$, $p\neq 1$, the additional
 operator
 \bb\label{e6} Y_5 = x\frac{\kd }{\kd x} +y\frac{\kd }{\kd
 y}+\frac{2}{1-p}u\frac{\kd }{\kd u} \ee
 generates a dilation.

\item For $f(u)=k e^u$, we have
 \bb\label{e7} Y_{(\xi^{1},\xi^{2})}^{e}={\xi }^1\frac{\kd }{\kd x} +{\xi }^2\frac{\kd }{\kd
 y} -2{\xi }^1_x\frac{\kd }{\kd u}, \ee
 where ${\xi }^1$ and ${\xi }^2$
 satisfy the Cauchy-Riemann system $(\ref{e8})$.
 \end{itemize}
 
Note that the projection of $Y_{(\xi^{1},\xi^{2})}^{e}$ on the $(x,y)$-space is the conformal group of $(\R^{2},ds^{2})$, where $ds^{2}=dx^{2}+dy^{2}$. For more details about two-dimensional Poisson equations, see \cite{yii}.

When $n>2$, the group classification is the same of the Polyharmonic equation taking $m=1$ in equation (\ref{pol}). See next section.

\subsection{Polyharmonic Equations}\label{poly}

The semilinear polyharmonic equation
\bb\label{pol}
(-1)^{m}\Delta^{m}u=f(u),
\ee
where $\Delta$ is the Laplace operator in $\R^{n},n\geq 2$ and $m\in\N$ is one of the most studied elliptic PDE. In \cite{sv}, Svirshchevskii proved that for any function $f(u)$, the widest Lie point symmetry group admitted by (\ref{pol}) is determined by translations and rotations, given, respectively, by the following vector fields in $\R^{n}$ 
\bb\label{ggenpol}
X_{i}=\f{\p}{\p x^{i}},\;\; X_{ij}=x^{j}\f{\p}{\p x^{i}}-x^{i}\f{\p}{\p x^{j}},\;\; 1\leq i,j\leq n,\;\;i\neq j.
\ee

In this paper, we consider equation (\ref{pol}) in $\R^{n}$ with $n>2$.

For special choices of function $f(u)$ in (\ref{pol}), the symmetry group can be enlarged. Below we exhibit these functions and their respective additional symmetries.

\begin{itemize}
\item If $f(u)=0$, the additional symmetries are 
\bb\label{ypol}Y_{i}=(2x^{i}x^{j}-\|x\|^{2}\delta^{ij})\f{\p}{\p x^{j}}+(2m-n)x^{i}u\f{\p}{\p u},
\ee
where $\delta^{ij}$ is the Kroenecker delta and $\|x\|$ is the Euclidean norm of $x$,
\bb\label{hompol}U=u\f{\p}{\p u},\;\;W_{\be}=\be\f{\p}{\p u},\ee
where $(-\Delta)^{m}\be=0$.

\item If $f(u)=u $, the additional summetries are $U$ and $W_{\be}$ as in (\ref{hompol}), and $\be$ satisfies $$(-1)^{m}\Delta^{m}\be+\be=0.$$ 

\item If $f(u)=u^{p} $, $p\neq 0, p\neq
 1,$ we have the generator of dilations
$$ D_{pm}= x^{i}\frac{\p}{\p x^{i}}+\f{2m}{1-p}u\frac{\p}{\p u}. $$
If $n\neq 2m$ and $p=(n+2m)/(n-2m)$, there are $n$ additional symmetries given by the vector fields (\ref{ypol}).

\item If $f(u)=e^{u}$ the additional symmetry is
$$W=x^{i}\f{\p}{\p x^{i}}-2m\f{\p}{\p u}$$
When $n=2m$, there are the following additional vector fields:
$$E_{i}=(2x^{i}x^{j}-\|x\|^{2}\delta^{ij})\f{\p}{\p x^{j}}-4m\f{\p}{\p u},$$
\end{itemize}
For more details about Group Analysis of equation (\ref{klgen}), see \cite{yb2, sv}.

\subsection{Wave Equations}
Hyperbolic type second-order nonlinear PDEs in two independent variables are used to describe different types of wave propagation. 

Consider the following semilinear wave equation in two independet variables
\bb\label{wavegen}
u_{tt}=u_{xx}+f(u).
\ee

For any function $f(u)$, the vector fields
\bb\label{simw}W_{1}=t\f{\p}{\p t}+x\f{\p}{\p x},\;\;W_{2}=\f{\p}{\p t},\;\; W_{3}=\f{\p}{\p x},\ee
are Lie point symmetries of equation (\ref{wavegen}). For some choices of functions $f(u)$, we have the following additional symmetries:
\begin{itemize}
\item If $f(u)=0$, then the symmetry group is
$$W_{\xi,\phi}=\xi(x,t)\f{\p}{\p x}+\phi(x,t)\f{\p}{\p t},\;\; U=u\f{\p}{\p u},\;\; \W=\beta\f{\p}{\p u},$$
where the functions $\xi,\phi,\beta$ satisfy
\bb\label{kilh}
\xi_{x}-\phi_{t}=0,\;\;
\xi_{t}-\phi_{x}=0,
\ee
$$
\be_{xx}-\be_{tt}=0.
$$ 
\item If $f=u$, then the symmetry group of (\ref{wavegen}) is generated by (\ref{simw}) and by
$$
U=u\f{\p}{\p u},\;\;W_{\beta}=\beta(x,t)\f{\p}{\p u},\;\;\text{where } \beta_{xx}-\beta_{tt}+\beta=0.
$$

\item If the nonlinearity is a power of $u$, i.e., $f(u)=u^{p}$, with $p\neq 0,1$, whe have the dilation symmetry
$$ D_{p}=x\f{\p}{\p x}+t\f{\p}{\p t}+\f{2}{1-p}u\f{\p}{\p u}.$$

\item If $f(u)=e^{u}$, then the symmetry group is
$$W_{\xi,\phi}^{e}=\xi\f{\p}{\p x}+\phi\f{\p}{\p t}-2\xi_{x}\f{\p}{\p u},$$
where $\xi,\phi$ satisfy (\ref{kilh}).
\end{itemize}

The projection of symmetry $W_{\xi,\phi}^{e}$ to the plane is the conformal group of $(\R^{2},ds^{2})$, where $ds^{2}=dx^{2}-dt^{2}$. It is analogous to the Euclidean case.

In \cite{la} there is a wide list of many kinds of wave equations. Here, we considered only a particular case. For more details, see \cite{la}.

\subsection{Heat Equations}
Consider the one-dimensional heat conduction equation
\bb\label{heat}
u_{t}=u_{xx}+f(u).
\ee
The symmetry group is generated by the following vector fields:
\begin{itemize}
\item For any function $f(u)$, the symmetries
\bb\label{heatgen}
H_{0}=\f{\p}{\p t},\;\;H_{1}=\f{\p}{\p x},
\ee
is a symmetru group of equation (\ref{heat}).
In addition to symmetries (\ref{heatgen}), for some choices of function $f(u)$, we have:
\item If $f(u)=0$, then
\bb\label{h2}
\ds{H_{2}=2t\f{\p}{\p x}-xu\f{\p}{\p u},\;\;H_{u}=u\f{\p}{\p u}},
\ee
$$
\ba{c}
\ds{H_{3}=x\f{\p}{\p x}+2t\f{\p}{\p t}},\\
\\
\ds{H_{4}=4tx\f{\p}{\p x}+4t^{2}\f{\p}{\p t}-(x^{2}+2t)u\f{\p}{\p u}},\\
\\
\ds{H_{\be}=\be\f{\p}{\p u}}, \text{ where }\ds{\be_{t}-\be{xx}=0}.
\ea
$$
\item If $f(u)=u$, we have the symmetries (\ref{h2}) and the following additional generators
$$
\ba{c}
\ds{H_{5}=x\f{\p}{\p x}+2t\f{\p}{\p t}+2tu\f{\p}{\p u}},\\
\\
\ds{H_{6}=tx\f{\p}{\p x}+t^{2}\f{\p}{\p t}+(t^{2}-\f{x^{2}}{4}-\f{t}{2})u\f{\p}{\p u}},\\
\\
\ds{H_{\be}=\be\f{\p}{\p u}}, \text{ where }\ds{\be_{t}-\be{xx}=0}.
\ea
$$
\item If $f(u)=u^{p},p\neq 0, 1,2$, we have
$$H^{d}_{p}=x\f{\p}{\p x}+2t\f{\p}{\p t}+\f{2}{1-p}u\f{\p}{\p u}.$$
\item If $f(u)=u^{2}$, then
$$
\ba{c}
\ds{H_{7}=tx\f{\p}{\p x}+t^{2}\f{\p}{\p t}-2tu\f{\p}{\p u}-\f{\p}{\p u}},\\
\\
\ds{H^{d}_{2}=x\f{\p}{\p x}+2t\f{\p}{\p t}-2u\f{\p}{\p u}}.
\ea
$$

\item If $f(u)=e^{u}$, then
$$H_{8}=x\f{\p}{\p x}+2t\f{\p}{\p t}-2\f{\p}{\p u}.$$

\end{itemize}

\subsection{Kohn - Laplace Equations}
The Heisenberg Group $H^{1}$ is the three-dimensional nilpotent Lie group, with composition law defined by  $$\R^{3}\times\R^{3}\ni((x,y,t),(x_{0},y_{0},t_{0}))\mapsto\phi((x,y,t),(x_{0},y_{0},t_{0})):=(x+x_{0}, y+y_{0}, t+t_{0}+2(xy_{0}-yx_{0}))\in\R^{3}.$$
In $H^{1}$ there is the subeliptic Laplacian defined by
$$\lh=X^2+Y^2,$$
where $\ds{X=\f{\p}{\p x}+2y\f{\p}{\p t}}$ and $\ds{Y=\f{\p}{\p y}-2x\f{\p}{\p t}}$.

The Kohn-Laplace equations is 
 \bb\label{klgen}
 u_{xx}+u_{yy}+4(x^{2}+y^{2})u_{tt}+4yu_{xt}-4xu_{yt}+f(u)=0,\ee
where $f:\R\rightarrow\R$ is a smooth function.

In \cite{gc} a complete group classification for equation (\ref{klgen}) is presented. It can be summarized as follows.

Let $G_{f}:=\{T, R,\X,\Y\}$, where 
$$ T=\frac{\p }{\p t},\;\;\; R= y\frac{\p }{\p
x}-x \frac{\p }{\p y},\;\;\; \tilde{X}=\frac{\p}{\p x}-2y \frac{\p}{\p
t}, \text{ and }\tilde{Y}=\frac{\p}{\p y}+2x \frac{\p}{\p t}.$$
For any function $f(u)$, the group $G_{f}$ is a symmetry group. 

For special choices of function $f(u)$ in (\ref{klgen}), the symmetry group can be enlarged. Below we exhibit these functions and their respective additional symmetries.

\begin{itemize}
\item If $f(u)=0$, the additional symmetries are
\bb\label{v1} V_1 = (xt-x^{2}y-y^{3})\frac{\p }{\p x} +
 (yt+x^{3}+xy^{2})\frac{\p }{\p y} + (t^{2}-(x^{2}+y^{2})^{2})\frac{\p }{\p t}-
 t u \frac{\p }{\p u},\ee
\bb\label{v2} V_2  = (t-4xy)\frac{\p }{\p x} +
 (3x^{2}-y^{2})\frac{\p }{\p y} - (2yt+2x^{3}+2xy^{2})\frac{\p }{\p t}
 + 2 y u \frac{\p }{\p u}, \ee
\bb\label{v3} V_3 = (x^{2}-3y^{2})\frac{\p }{\p x} +
 (t+4xy)\frac{\p }{\p y} + (2xt-2x^{2}y-2y^{3})\frac{\p }{\p t}
 - 2x u \frac{\p }{\p u},\ee
$$ Z =x\frac{\p}{\p x}+y\frac{\p}{\p y}+2t
\frac{\p}{\p t} ,  \;\;\; U= u\frac{\p }{\p u},\;\;\;\W= \be
(x,y,t)\frac{\p }{\p u}, \text{ where } \lh \be =0.
$$

\item If $f(u)=u $, there are two additional symmetries 
$$U = u\frac{\p }{\p u},\;\;\; \W= \be (x,y,t)\frac{\p }{\p u}, \text{  where } \lh \be +\be=0.$$

\item If $f(u)=u^{p} $, $p\neq 0, p\neq
 1, p\neq 3$, we have the generator of dilations
\bb\label{d3} D_{p}= x\frac{\p}{\p x}+y\frac{\p}{\p y}+2t
 \frac{\p}{\p t}+\frac{2}{1-p}u\frac{\p}{\p u}. \ee

\item If $f(u)=e^{u}$ the additional symmetry is
$$ E=x\frac{\p}{\p x}+y\frac{\p}{\p y}+2t
\frac{\p}{\p t}-2\frac{\p}{\p u}. $$ 

\item In the critical
 case, $f(u)=u^{ 3} $, there are four additional generators, namely
 $V_1,V_2,V_3\text{ and } D_{3}$, given in $(\ref{v1}), (\ref{v2}),(\ref{v3})\text{ and }(\ref{d3})$
 respectively. 
\end{itemize}
For more details about Group Analysis of equation (\ref{klgen}), see \cite{gc, ds, cl, ilf}.





\begin{center}\textbf{Acknowledgements}\end{center}

\

I. L. Freire would like to thank I. I. Onnis, W. A. Rodrigues Jr and M. A. Faria Rosa for their encouragement to prove Theorem 1. He is grateful to UNICAMP for financial support. We also thank Lab. EPIFISMA (Proj. FAPESP) for having given us
the opportunity to use excellent computer facilities.  

\end{document}